\documentclass[twoside,12pt]{amsart}

\usepackage{enumerate}
\usepackage{graphics}
\usepackage{graphicx}
\usepackage{amsmath}
\usepackage{amssymb}
\usepackage{amstext}
\usepackage{times}
\usepackage[latin1]{inputenc}
\usepackage{latexsym}
\usepackage{amsthm}
\usepackage{amscd}

\def\R{\mathbb{R}}
\def\C{\mathbb{C}}

\def\ep{{\epsilon}}

\def\a{{\alpha}}
\def\g{{\gamma}}
\def\be{{\beta}}
\def\l{{\lambda}}
\def\grad{{\nabla}}
\def\ri{{\rm i}}

\newtheorem{lemma}{Lemma}
\newtheorem{remark}{Remark}
\newtheorem{theorem}{Theorem}

\newtheorem{definition}{Definition}

\numberwithin{equation}{section}

\begin{document}

\title{Limit sets for complete minimal immersions}
\author{Antonio Alarcón}
\author{Nikolai Nadirashvili}

\thanks{A. Alarcón's research is partially supported by MEC-FEDER Grant no. MTM2004 - 00160.}

\date{\today}

\address{Antonio Alarcón, \hfill\break\indent Departamento de Geometría y Topología,\hfill\break\indent  Universidad de Granada, \hfill\break\indent 18071, Granada, \hfill\break\indent Spain.}
\email{alarcon@ugr.es}

\address{Nikolai Nadirashvili, \hfill\break\indent CNRS, LATP, CMI,\hfill\break\indent 39, rue Joliot-Curie, \hfill\break\indent 13453 Marseille Cedex 13, \hfill\break\indent France.}
\email{nicolas@cmi.univ-mrs.fr}

\begin{abstract} 
In this paper we study the behaviour of the limit set of complete proper compact minimal immersions in a domain $G\subset \R^3$ with the boundary $\partial G\subset C^2.$ We prove that the second fundamental form of the surface $\partial G$ is nonnegatively defined at every point of the limit set of such immersions.
\\


\end{abstract}

\maketitle


\section{Introduction}\label{sec: intro}

Let $F:M\to\R^3$ be a complete conformal minimal immersion. After the discovery in 1996 by the second author \cite{Na1} that the surface can be bounded, i. e., the image of $F$ can be in a ball of $\R^3,$ it was done a lot of work studying the topology and limit sets of bounded complete minimal immersions (see for instance \cite{LMM1}, \cite{MM2} and \cite{AFM}). In \cite{MM1} it was proven that as a surface $M$ can be taken a disk. From now on, we will assume that $M$ is the unite disk of $\C.$

Denote by $\g$ the limit set of the surface $F(M),$ i. e., $\g\subset \R^3$ is the limit set of $F(z)$ for $z\to \partial M.$ It was proven \cite{MM3} that for any bounded convex domain $G\subset \R^3$ there exists a complete minimal surface such that its limit set $\g\subset \partial G.$ On the other hand, in \cite{MMN} was shown the existence of a domain $G\subset \R^3$ for which there is no complete properly immersed in $G$ minimal surfaces.

If we assume that the minimal surface is in addition compact, i. e., the map $F$ has a continuous extension to the map $F:\overline{M}\to\R^3,$ then $\g$ can not be a subset of the boundary of a cube in $\R^3$ \cite{Na2}. Unfortunately in the paper \cite{Na2} the condition of compactness of the minimal surface was missed in the statement of the theorem, however was used in the proof.
Notice that for compact minimal immersions the limit set $\g$ coincides with the set $F(\partial M).$ In \cite{MN} it was proven that there are compact complete bounded minimal immersions such that $\g$ is a Jordan curve of Hausdorff dimension 1. It is easy to see that for complete minimal surfaces $\g$ can not be a rectifiable curve.
\\

In this paper we prove the following results:

\begin{theorem}\label{th: A}
Let $G\subset\R^3$ be a domain with the boundary $\partial G\subset C^2.$ Let $F:M\to G$ be a proper compact complete minimal immersion. Then $F(\partial M)\subset \partial G$ and if $x\in F(\partial M)$ then the second fundamental form of the surface $\partial G$ at $x$ is nonnegatively defined.
\end{theorem}

\begin{theorem}\label{th: B}
There is no properly immersed compact complete minimal disk in a polyhedron $P\subset \R^3$ nor in $\R^3\setminus \overline{P}.$
\end{theorem}


\section{Preliminaries}\label{se: pre}

As a previous step to prove the above results we are going to introduce some notation.

\begin{definition}\label{def: concave}
Consider $G\subset\R^3$ a domain with the boundary $\partial G\subset C^2,$ and a point $p \in \partial G.$ Then, we say that $p$ is a concave point of $G$ if the principal curvatures of $\partial G$ associated to the inward pointing unit normal are not positive on a neighbourhood of $p.$
\end{definition}

From this definition, we obtain the following two results straightforwardly:

\begin{remark}\label{rem: entorno}
Let $G\subset\R^3$ be a domain with the boundary $\partial G\subset C^2,$ and $p\in \partial G$ a concave point of $G.$ Then, there exists a neighbourhood $U$ of $p$ in $\partial G$ so that $q$ is a concave point of $G$ for all $q\in U.$
\end{remark}

\begin{remark}\label{rem: plane}
Notice that given a such domain $G$ of $\R^3$ and a concave point $p$ of $G$, then the tangent plane $T_p (\partial G)$ satisfies $p \in T_p (\partial G)$ and there exists a neighbourhood of $p$ in $T_p (\partial G)$ contained in $\overline{G}.$
\end{remark}

In order to prove the main theorems, besides these two immediate results, we will use a deep result of Bourgain \cite{Bo}.

\begin{theorem}[Bourgain]\label{th: Bourgain}
Consider $u:M\to\R$ a bounded harmonic function. Then, there exists an open dense set $A$ on $[0,2\pi[$ so that the integral
\[
\int_0^1 |\grad u(r e^{\ri \a})|\, dr
\]
is convergent for all $\a\in A.$
\end{theorem}

Moreover, during the proof of Theorem \ref{th: A}, we will need the following two technical lemmas:

\begin{lemma}\label{lem: plane}
Let $G\subset\R^3$ be a domain with the boundary $\partial G\subset C^2,$ and $p\in \partial G$ a concave point of $G.$ Let $X:M\to G$ be a compact proper minimal immersion (not necessarily complete) so that $p=X(e^{\ri \be}),$ where $\be\in\R.$ Label $\Pi$ the tangent plane to $G$ at $p.$ Then, there exists an open interval $I\subset\R$ satisfying $\be\in I$ and such that given $\a_1,\a_2\in I$ with $\a_1<\be<\a_2$ and a curve $\gamma\subset M$ joining $e^{\ri \a_1}$ and $e^{\ri \a_2},$  then $X(\g)\cap\Pi\neq \emptyset.$
\end{lemma}

\begin{lemma}\label{lem: preimage}
Let $u:M\to \R$ be a harmonic function having a continuous extension to $u:\overline{M}\to\R.$ Consider $\be \in \R$ satisfying the following two properties:
\begin{itemize}
\item There exists an open interval $I\subset\R$ such that $\be \in I$ and $u_{|\sigma}\geq 0,$ where $\sigma = \{e^{\ri \a}\;|\; \a\in I\}.$
\item There exists an arc $\g\subset\overline{M}$ such that $e^{\ri \be}\in\g$ and $u_{|\g}\equiv 0.$
\end{itemize} 
Then,
\[
\int_0^1 |\grad u(r e^{\ri \be})|\, dr<\infty\;.
\]
\end{lemma}

Lemmas \ref{lem: plane} and \ref{lem: preimage} will be proved in Section \ref{sec: lemma}.


\section{Proof of Theorem \ref{th: A}}\label{sec: proofA}

Theorem \ref{th: A} is a trivial consequence of the following one:

\begin{theorem}\label{th: regular}
Let $G\subset\R^3$ be a domain with the boundary $\partial G\subset C^2$ and $X:M\to G$ a proper compact  complete minimal immersion. Then, the limit set of $X$ does not contain concave points of $G.$
\end{theorem}
\begin{proof}
We will suppose that there exists a concave point of $G$ in the limit set of $X$ and we will lead us to a contradiction.

Assume the existence of such a point $p.$ Since $p$ belongs to the limit set, there exists $\a \in \R$ so that $p=X(e^{\ri \a}).$ $X$ is proper and compact, hence, since Remark \ref{rem: entorno} we can find an open interval $I\subset\R$ so that $\a\in I$ and $X(e^{\ri \be})$ is a concave point of $G$ for all $\be\in I.$ Up to a rigid motion, we can assume that
\begin{equation}\label{eq: tangent}
T_p(\partial G)=\{(x_1,x_2,x_3)\in\R^3\;|\;x_1=1\}\;.
\end{equation}

Using again that $X\equiv (x_1,x_2,x_3)$ is compact, we obtain that it is bounded. Therefore, $x_2:M\to\R$ is a bounded harmonic function. Hence, the result of Bourgain (Theorem \ref{th: Bourgain}) guarantees the existence of a dense set $A$ on $I$ so that
\begin{equation}\label{eq: x2}
\int_0^1 |\grad x_2 (r e^{\ri \be})|\, dr<\infty\;, \quad \forall \be\in A\;.
\end{equation}
Choose $\be\in A$ close enough to $\a$ so that
\[
T_q(\partial C)=\{(x_1,x_2,x_3)\in\R^3\;|\; \sum_{i=1}^3\l_ix_i=1\}\;,
\]
where $\l_1\neq 0,\l_2,\l_3\in \R$ and we have denoted $q=X(e^{\ri \be}).$ This election is possible because of \eqref{eq: tangent}. Now, consider the function
\[
v=\sum_{i=1}^3 \l_ix_i:\overline{M}\to \R\;,
\]
which is continuous on $\overline{M},$ harmonic on $M$ and identically equal to $1$ on the set $X^{-1}(T_q(\partial G)).$ At this point notice that Lemma \ref{lem: plane} guarantees that $X^{-1}(T_q(\partial G))$ is an arc. Moreover, $e^{\ri\be}\in X^{-1}(T_q(\partial G)).$ 
Therefore, taking Remark \ref{rem: plane} into account, we can make use of Lemma \ref{lem: preimage} to obtain that the integral
\[
\int_0^1 |\grad v (r e^{\ri \be})|\, dr < \infty\;.
\]
\begin{figure}[htbp]
	\begin{center}
		\includegraphics[width=0.80\textwidth]{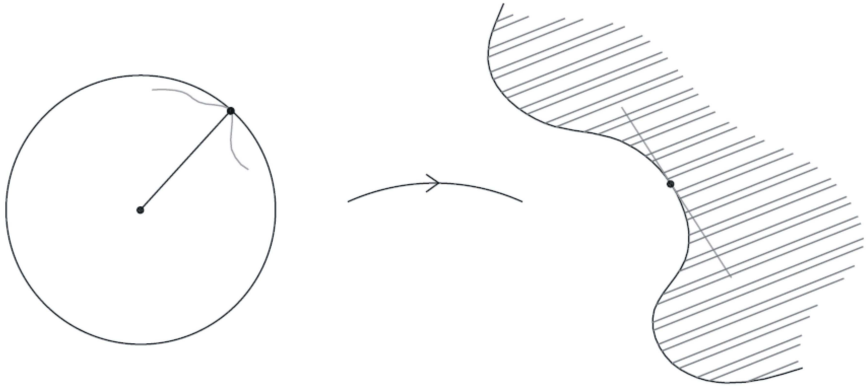}
	\end{center}
	\caption{The curve $X^{-1}(T_q(\partial G)).$} \label{fig: curve}
\end{figure}

Consequently, using also \eqref{eq: x2} and labeling $\l=\l_3/\l_1,$ we have
\begin{equation}\label{eq: v}
\int_0^1 |\grad x_1(r e^{\ri \be})+ \l \grad x_3(r e^{\ri \be})|\, dr < \infty\;.
\end{equation}

On the other hand, $X$ is a minimal immersion of a Riemannian surface, and so, it is a conformal map. From this fact it follows that $(dx_1)^2+(dx_2)^2+(dx_3)^2$ is a conformal metric on $M.$ Then, working with $\grad x_i$ as complex numbers, we have
\begin{equation}\label{eq: conformal}
(\grad x_1)^2+(\grad x_2)^2+(\grad x_3)^2 = 0\;.
\end{equation}
\begin{figure}[htbp]
	\begin{center}
		\includegraphics[width=0.50\textwidth]{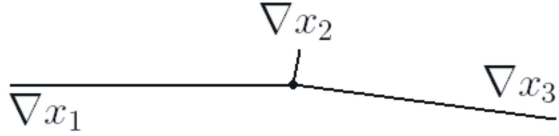}
	\end{center}
	\caption{This kind of situation is not a conformal map.} \label{fig: conformal}
\end{figure}

Now, taking into account \eqref{eq: x2}, \eqref{eq: v} and \eqref{eq: conformal}, we obtain that
\[
\int_0^1 |\grad x_i(r e^{\ri \be})|\; dr<\infty\;,\quad \forall i=1,2,3\;.
\]
As a consequence, we have
\[
\int_0^1 |\grad X (r e^{\ri \be})|\, dr<\infty\;,
\]
that means that the curve $X(r e^{\ri\be}),$ $0<r<1,$ has finite length, and hence the minimal surface $X(M)$ is not complete. This contradicts the hypotheses of the theorem, and so, finishes the proof.
\end{proof}




\section{Proof of Theorem \ref{th: B}}\label{sec: proofB}

The proof is quite similar to the proof of Theorem \ref{th: regular}. Therefore, we will not explain some of the details. Again, we proceed by contradiction. Assume there exists a such immersion $X=(x_1,x_2,x_3).$ Up to a rigid motion, the compactness of $X$ guarantees the existence of an open interval $I\subset\R$ so that
\[
x_1(e^{\ri \a})=1\;, \quad \forall \a\in I\;.
\]
Hence, $x_1$ is a smooth function in a neighbourhood of $I,$ and so
\begin{equation}\label{eq: x1'}
\int_0^1 |\grad x_1 (r e^{\ri \a})|\, dr<\infty\;, \quad \forall \a\in I\;.
\end{equation}

On the other hand, since $x_2$ is a bounded harmonic function, Theorem \ref{th: Bourgain} gives us $\be\in I$ so that
\begin{equation}\label{eq: x2'}
\int_0^1 |\grad x_2 (r e^{\ri \be})|\, dr<\infty\;.
\end{equation}

On the other side, $X$ is a conformal map, thus
$$
|\grad x_3| \leq |\grad x_1| + |\grad x_2|\;.
$$
and so, taking into account \eqref{eq: x1'} and \eqref{eq: x2'}, we conclude that
\begin{equation}\label{eq: x3'}
\int_0^1 |\grad x_3 (r e^{\ri \be})|\, dr<\infty\;.
\end{equation}

Finally, combining the inequalities \eqref{eq: x1'}, \eqref{eq: x2'} and \eqref{eq: x3'} we obtain
\[
\int_0^1 |\grad X (r e^{\ri \be})|\, dr<\infty\;,
\]
which proves the theorem.


\section{Proof of the technical lemmas}\label{sec: lemma}

\subsection{Proof of Lemma \ref{lem: plane}}
We will proceed by contradiction. Assume the existence of real numbers $\a_1$ and $\a_2$ with $\a_1<\be<\a_2$ and a curve $\gamma\subset M$ joining $e^{\ri\a_1}$ and $e^{\ri\a_2}$ and such that $X(\g)\cap \Pi= \emptyset.$ Label $\gamma'=\{e^{\ri \delta}\;|\;\a_1\leq \delta\leq \a_2\}$ and let $D$ be the open subset of $M$ bounded by $\gamma\cup\gamma'.$ Then, the convex hull property of minimal surfaces (see \cite{Os2}) guarantees that $X(D)\subset E,$ where $E$ is convex hull of the set $\gamma\cup\gamma'.$ Hence, from our hypothesis, we have
\begin{equation}\label{eq: subset}
X(D)\subset G\cap E\;.
\end{equation}

Now, notice that if $\a_1$ and $\a_2$ are close enough to $\be,$ then, since Remark \ref{rem: entorno}, $\gamma'$ is a curve of concave points of $\partial G.$ Since this fact and taking into account that $X(\g)\cap \Pi= \emptyset,$ we obtain the existence of a small neighbourhood $E'$ of $p$ in $E$ so that $G\cap E'=\emptyset.$ This contradicts \eqref{eq: subset} and proves the lemma.

\subsection{Proof of Lemma \ref{lem: preimage}}
First of all, notice that if a neighbourhood of $e^{\ri \be}$ in $\gamma$ lies on $\partial M,$ then $u$ is a differentiable function in a neighbourhood of $e^{\ri \be}$ and so, the lemma trivially holds. Therefore, without loss of generality, we assume $\gamma\setminus\{ e^{\ri\be}\}\subset M.$
\\

Now, define the functions
\[
\varphi_1 =\sup \{0,u_{|\partial M}\}\;,\quad\text{and}\quad \varphi_2=\inf\{ 0,u_{|\partial M}\}\;.
\]
For each $i=1,2,$ consider $v_i$ the function solving the problem
\[
\begin{cases}
\Delta v_i=0 & \text{in } M \\
v_i=\varphi_i & \text{in } \partial M\;.
\end{cases}
\]
Then, we have
\begin{equation}\label{eq: sum}
u=v_1+v_2\;,
\end{equation}
and $v_2\leq 0 \leq v_1$ in $M.$ Moreover, since ${v_2}_{|\sigma}=0$ we obtain that fixed $\ep>0$ there exists a constant $C>0$ so that
\begin{equation}\label{eq: bound}
-v_2(x) < C |x-e^{\ri\be}|\quad \forall x\in M \text{ with } |x-e^{\ri\be}|< \ep\;.
\end{equation}
Therefore, since $v_2$ is harmonic,
\begin{equation}\label{eq: v2}
|\grad v_2(x)| < C \quad \forall x\in M \text{ with } |x-e^{\ri\be}|< \ep\;.
\end{equation}

Now, for $0<\rho<1/2$ define $x_\rho = (1-2\rho)e^{\ri \be}.$ Let $B_\rho$ be the ball centered at $x_\rho$ and with radius $\rho.$
\begin{figure}[htbp]
	\begin{center}
		\includegraphics[width=0.60\textwidth]{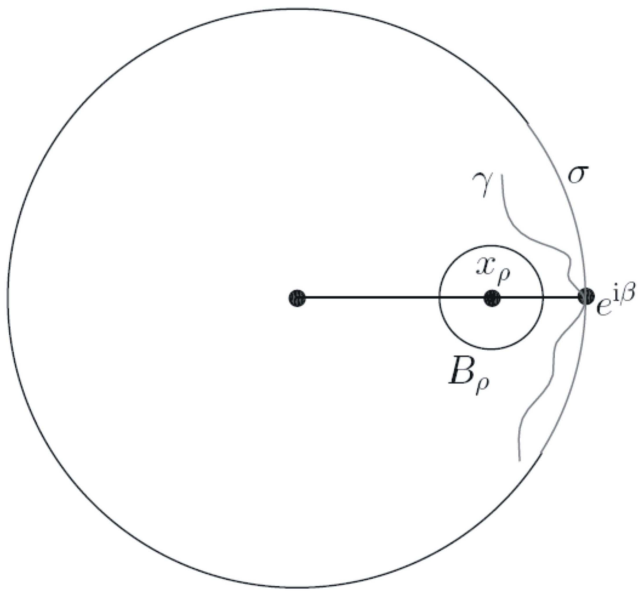}
	\end{center}
	\caption{The curves of Lemma \ref{lem: preimage} and the ball $B_\rho.$} \label{fig: preimage}
\end{figure}

Fix $C_0>0$ and assume that $C_0 v_1(x_\rho) > \rho.$ Then, since $v_1\geq 0$ in $M,$ by Harnack's Theorem we have a constant $C_1>0$ (which does not depend on $C_0$) so that
\[
C_0\, {v_1}_{|B_\rho}> C_1 \rho\;.
\]

On the other hand, since $v_1\geq 0$ in $\partial M,$ we know that there exists another constant $C_2>0$ (again, it does not depend on $C_0$) satisfying
\[
C_0\, v_1(x)> C_2\, G(x_\rho,x)\quad \forall x \in \partial B_\rho\;,
\]
where $G$ is the Green function in $M,$ having its singularity in $x_\rho.$ Hence,
\[
C_0\, v_1(x) > C_2\, G(x_\rho,x)\quad \forall x\in M\setminus B_\rho\;.
\]
This fact implies the existence of a new constant $C_3>0$ (not depending on $C_0$) such that
\[
C_0\, v_1(x) > C_3\, |x-e^{\ri\be}|\;.
\]
Therefore, taking $C_0$ small enough so that $C_3>C_0 C,$ where $C$ is the constant of \eqref{eq: bound}, we obtain a contradiction with the fact that $u_{|\g}=0.$ Thus, there exists a constant $C_4>0$ such that $v_1< C_4 \rho.$ Hence, there exists $C'>0$ satisfying
\begin{equation}\label{eq: v1}
|\grad v_1(x_\rho)|< C'\;.
\end{equation}

Finally, we conclude the proof combining \eqref{eq: sum}, \eqref{eq: v1} and \eqref{eq: v2}.


\end{document}